\documentclass[12pt]{article}

\def\be{\begin{equation}}
\def\ee{\end{equation}}
\def\bea{\begin{eqnarray}}
\def\eea{\end{eqnarray}}
\def\bean{\begin{eqnarray*}}
\def\eean{\end{eqnarray*}}

\def\la{\label}
\def\r{\ref}

\begin{document}
\rightline{June 9, 2005}
\vskip 1 true cm
\begin{center}
{\Large Gaussian Quadrature without Orthogonal Polynomials}
\vskip 0.6 true cm
{\large Ilan Degani}\\
\smallskip
{\it Department of Chemical Physics}\\
{\it Weizmann Institute of Science, Rehovot 76100, Israel}\\
{\small ilan.degani@weizmann.ac.il}\\
{\small\em and}\\
{\large Jeremy Schiff}\\
\smallskip
{\it Department of Mathematics}\\
{\it Bar-Ilan University, Ramat Gan 52900, Israel}\\
{\small schiff@math.biu.ac.il}\\
\end{center}
\vskip 0.6 true cm

\begin{quote}
\centerline{\bf Abstract}

A novel development is given of the theory of 
Gaussian quadrature, not relying on the 
theory of orthogonal polynomials. A method
is given for computing the nodes and weights that is 
manifestly independent of choice of basis in the space of polynomials.
This method can be extended to compute nodes and weights for 
Gaussian quadrature on the unit circle and Gauss type quadrature 
rules with some fixed nodes. 

\end{quote}

\vskip 0.6 true cm

The aim of this letter is to show how the theory of 
Gaussian quadrature can be developed using elementary 
linear algebra, without relying on the theory of 
orthogonal polynomials. This is certainly useful from a 
pedagogical point of view, showing how basic results in 
linear algebra can have powerful applications and
allowing Gaussian quadrature to be taught in introductory 
courses before the rather heavier topic of
orthogonal polynomials. But it is also 
important from a more fundamental point of view. Orthogonal 
polynomials are an invaluable technical tool, but they are,
after all, just a  choice of basis 
in the vector space of polynomials. 
In this paper we give a method
for computing nodes and weights that is 
independent of the choice of basis.

A number of the formulae in this letter can also be found in the papers
of Sack and Donovan [7] and Gautschi [5] on the ``modified moment 
method'' for computing nodes and weights. The emphasis, however,
is somewhat different: these papers focus on efficient, 
stable computation, while the focus here is more on the conceptual 
issue of basis indepedence. It also turns out that occasional benefit 
can be obtained by a choice of basis not covered in [7] and [5]. 
Our basis independent method for computing nodes and weights is described
in theorem 3. It allows immediate generalization
to the cases of Gaussian quadrature on the unit circle and 
Gauss type quadrature with some fixed nodes, as we shall show.
Our view of quadrature formulae as being related to blinear forms 
allows a generalization to multiple dimensions, though this is
discussed in a separate paper [3].

There is, of course, already a completely elementary approach to Gaussian 
quadrature that avoids mention of orthogonal polynomials.
We seek $x_i,w_i$, $i=1,\ldots,n+1$, such that
the quadrature formula
\be
\int_a^b w(x) f(x) dx \approx \sum_{i=1}^{n+1} w_if(x_i) \ , \qquad w_i\not=0\ ,
\label{quad}
\ee
is exact when $f(x)$ is any polynomial of degree not more than 
$2n+1$. Here $w(x)$ is a known non-negative weight function on 
the interval $[a,b]$.
Substituting in turn $f(x)=1,x,x^2,\ldots,x^{2n+1}$ we 
see that the nodes and weights must satisfy the equations 
\be
\sum_{i=1}^{n+1} w_i x_i^{\alpha} = m_\alpha\ , \qquad \alpha=0,1,\ldots,2n+1\ ,
\la{direqs}\ee
where $m_\alpha$ is the $\alpha$th moment
\be
m_\alpha = \int_a^b w(x) x^\alpha  dx \ .
\la{mom}\ee
In principle we should be able to obtain the nodes and weights by 
solving the system (\r{direqs}) of nonlinear equations. 
The drawback of this approach is that there is no evident way to solve
(\r{direqs}) except in very simple  cases;
it is far from clear that there even exists a solution. Also, dependence 
on choice of basis is hardly resolved in this approach.

\smallskip

\noindent{\bf A linear algebra  approach to Gaussian quadrature.} Let
${\cal P}_n$ denote the vector space of polynomials of degree up to $n$
(so ${\rm dim}{\cal P}_n=n+1$). Let $w(x)$ be a given non-negative 
weight function on the interval $[a,b]$, and write
\be 
<f|g> = \int_a^b w(x) f(x) g(x)\ dx  \la{ip}
\ee
for any polynomials $f,g$. For any $n$, $<\cdot|\cdot>$ determines an 
inner product on the space ${\cal P}_n$. Define the symmetric bilinear form 
$X(\cdot,\cdot)$ on ${\cal P}_n$ by
\be
 X\left(f,g\right)  = \int_a^bw(x)x f(x) g(x)\  dx\ . \la{X}
\ee
For any symmetric bilinear form on an inner product space there
exists an orthonormal basis in which the form is diagonal. In other words, we can
find eigenvalues $x_1,\ldots,x_{n+1}\in{\bf R}$ 
 and eigenvectors $u_1(x),\ldots,u_{n+1}(x)\in{\cal P}_n$
such that 
\bea
< u_i | u_j > &=& \delta_{ij}\ , \\
X( u_i, u_j) &=& x_i \delta_{ij}\ .
\eea

\noindent {\bf Theorem 1}. 
Let $x_1,\ldots,x_{n+1}$ be the eigenvalues  of the bilinear form $X(\cdot,\cdot)$ with
corresponding orthonormal eigenvectors $u_1(x),\ldots,  u_{n+1}(x)$. Then for 
any polynomial $p(x)\in{\cal P}_{2n+1}$
\be
\int_a^b w(x)p(x) dx = \sum_{i=1}^{n+1} <1|u_i>^2 p(x_i)\ ,
\ee
i.e. the quadrature rule 
\be
\int_a^b w(x)f(x) dx \approx \sum_{i=1}^{n+1} w_i f(x_i)\ ,
\ee
where $w_i=<1|u_i>^2$, is exact whenever $f(x)\in{\cal P}_{2n+1}$. 

\noindent The proof proceeds via a lemma. 

\noindent {\bf The $\delta$ Lemma}. For any polynomial $p\in{\cal P}_{n+1}$,
we have
\be <p|u_i>=<1|u_i>p(x_i)\ .\la{dl}\ee
In other words, up to a normalization factor,
the projection of $p(x)$ on $u_i(x)$ is found by simply evaluating $p(x)$ at
$x=x_i$. 

\noindent {\bf Proof of The $\delta$ Lemma}. We prove the result for 
$p(x)=1,x,\ldots,x^{n+1}$,
the general result follows by linearity. For $p(x)=1$ the result is trivial. For
$j=1,\ldots,n+1$ we have 
\be
<x^j|u_i> = \int_a^b w(x) x^j u_i(x) \ dx 
             = X( x^{j-1}, u_i )  = x_i <x^{j-1}|u_i> \ .
\ee
Iterating this procedure $j$ times we obtain $<x^j|u_i>=x_i^j<1|u_i>$, 
as required. $\bullet$

\noindent{\bf Proof of Theorem 1}. We prove the result for $p(x)=1,x,\ldots,x^{2n+1}$,
the general result follows by linearity. For $p(x)=1$ we have 
\be
\int_a^b w(x)\ dx = <1|1> = \sum_{i=1}^{n+1} <1|u_i><u_i|1> =
    \sum_{i=1}^{n+1} <1|u_i>^2 
\ee
as required. 
For $p(x)=x^j$, $j=1,\ldots,2n+1$,  we can write  $p(x)=x^{j_1+j_2+1}$, 
where $j_1,j_2\in\{1,2,\ldots,n\}$, and thus
\bea
\int_a^b w(x)x^j \ dx &=&  X(x^{j_1},x^{j_2}) \nonumber \\ 
  &=& 
\sum_{i_1=1}^{n+1} \sum_{i_2=1}^{n+1} <x^{j_1}|u_{i_1}> <x^{j_2}|u_{i_2}> 
 X\left(u_{i_1},u_{i_2}\right)  \nonumber \\
&=&    \sum_{i=1}^{n+1} <1|u_i>^2 x_i^{j} \ ,
\eea
where in the last step we have used the $\delta$ Lemma. $\bullet$

Theorem 1 proves the existence of a Gaussian quadrature formula for any 
interval $[a,b]$ and non-negative weight function $w(x)$. The proof is based
just on diagonalization of the bilinear form $X(\cdot,\cdot)$ with no prior
knowledge needed of orthogonal polynomials. The next
theorem shows that some of the main results on Gaussian quadrature are easily
obtained within our approach.

\noindent{\bf Theorem 2}.
\newline 1) The Gaussian quadrature formula of theorem 1 is unique. 
\newline 2) The nodes of the Gaussian quadrature formula lie in the interval $[a,b]$.
\newline 3) The nodes of the Gaussian quadrature formula are roots of any 
nontrivial degree $n+1$ polynomial that is orthogonal to all polynomials 
of degree at most $n$. 

\noindent To motivate the first part of the proof
we make the observation that by setting 
$p(x)=u_j(x)$ in the $\delta$ Lemma we obtain
\be 
u_j(x_i) = \frac{\delta_{ij}}{<1|u_i>}\ ,
\ee
implying the interesting result
that the functions $u_i(x)$ are actually proportional to the 
Lagrange cardinal functions for the points $x_1,\ldots,x_{n+1}$. 
Note that by setting $p(x)=u_i(x)$ in (\r{dl}) we see that 
$<1|u_i>$ cannot be zero.
Note also that the last equation implies 
the $x_i$ are distinct. We mention in passing that 
the fact that $u_i(x_j)$ is proportional to $\delta_{ij}$ is one reason
we gave the $\delta$ Lemma its name.

\noindent{\bf Proof}.
1) Suppose $X_i,W_i$, $i=1,\ldots,n+1$ are nodes and weights for a 
Gaussian quadrature formula, i.e. for any polynomial $p(x)\in{\cal P}_{2n+1}$
\be
\int_a^b w(x)p(x)\ dx = \sum_{i=1}^{n+1} W_i\  p(X_i)\ ,
\ee
with the $X_i$ distinct and $W_i>0$.  
Define the $n$th degree polynomials $U_1(x),\ldots,U_{n+1}(x)$ by 
requiring
\be 
U_j(X_i) = \frac{\delta_{ij}}{\sqrt{W_i}}\ , \qquad i,j=1,\ldots,n+1 \ .
\ee
We then have 
\be
<U_i|U_j> = 
\int_a^b w(x) U_i(x)U_j(x) \ dx 
= \sum_{k=1}^{n+1} W_k U_i(X_k) U_j(X_k) = \delta_{ij}
\ee
and
\be
X\left(U_i|U_j\right) = 
\int_a^b w(x) x U_i(x)U_j(x) \ dx 
= \sum_{k=1}^{n+1} W_k X_k U_i(X_k) U_j(X_k) = X_i\delta_{ij}\ .
\ee
Thus the $X_i$ are necesarilly 
the eigenvalues of the bilinear form $X$, and the $U_i$ are
the orthonormal eigenvectors. Furthermore we have
\be
<1|U_i>= 
\int_a^b w(x) U_i(x) \ dx 
= \sum_{k=1}^{n+1} W_k U_i(X_k) = \sqrt{W_i} \ ,
\ee
so $W_i=<1|U_i>^2$, as before. 
\newline 2) We have 
$$ \int_a^b w(x) a u_i^2(x)  \ dx\le
    \int_a^b w(x) x u_i^2(x) \ dx \le
    \int_a^b w(x) b u_i^2(x) \ dx\ , $$
so $a\le X(u_i,u_i)\le b $ and therefore $a\le x_i\le b$.  
\newline 3) Note that the $\delta$ Lemma holds for all $p(x)\in{\cal P}_{n+1}$. 
If $p(x)$ is a polynomial of degree $n+1$ that is orthogonal to all polynomials
of degree at most $n$, then it is orthogonal to all the $u_i(x)$, and thus,
by the $\delta$ Lemma, $p(x_i)=0$
for all $i$. 
$\bullet$

\noindent{\bf Calculation of the nodes and weights.}
Having presented our nonstandard development of 
the theory of Gaussian quadrature, we turn our attention 
to the question of 
computing nodes and weights using an algorithm that allows
for an arbitrary choice of basis in ${\cal P}_n$.
To find the nodes and weights we clearly need information
on the inner product $<\cdot|\cdot>$ and the blinear form $X(\cdot,\cdot)$
on ${\cal P}_n$. So assume the functions $q_1(x),q_2(x),\ldots,q_{n+1}(x)$
are a basis of ${\cal P}_n$ and that we are given the 
$(n+1)\times(n+1)$ matrices $B,A$  defined by
\bea
B_{ij}&=&<q_i|q_j>=\int_a^b w(x) q_i(x) q_j(x) dx\ ,\qquad i,j=1,\ldots,n+1\ , \\
A_{ij}&=&X(q_i,q_j)=\int_a^b w(x) x q_i(x) q_j(x) dx\ ,\qquad i,j=1,\ldots,n+1\ .
\eea
We are about to see that knowledge of the matrices $B,A$ and just
one of the functions $q_i(x)$ is sufficient to determine the nodes
$x_i$ and weights $w_i$. As a prelude, note that there exists
a diagonal matrix $D$ and a matrix $V$ satisfying
\be 
V^TBV = I \ , \qquad  A V = B V D \ . \la{VD}
\ee
This is nothing but the statement that $X(\cdot,\cdot)$ has orthonormal
eigenvectors, written in the basis $\{q_i\}$. 
In Matlab Release 14 the matrices $D,V$ are computed from $A,B$
by the command {\tt [V D]=eig(A,B,'chol')}.

\noindent{\bf Theorem 3}.  If $D,V$ are matrices satisfying (\r{VD}), 
with $D$ diagonal, then
the nodes and weights of the Gaussian quadrature formula (that is exact
for polynomials in ${\cal P}_{2n+1}$) are determined by
\be
x_i = D_{ii}\ , \qquad
w_i = \left( \frac{(V^{-1})_{ij}}{q_j(x_i)} \right)^2 \ .
\ee
(The second relation holds for arbitrary $j$). 

\noindent{\bf Proof.}
Since all the integrals needed to compute $B$ and $A$ are integrals
of polynomials of degree up to $2n+1$, we have 
\bea
B_{ij}&=&\sum_{k=1}^{n+1} w_k q_i(x_k) q_j(x_k) = 
    \sum_{k=1}^{n+1} Q_{ik}W_{kk} Q_{jk} \ ,\la{BQ}\\
A_{ij}&=&\sum_{k=1}^{n+1} w_k x_k q_i(x_k) q_j(x_k) = 
    \sum_{k=1}^{n+1} Q_{ik} W_{kk}X_{kk} Q_{jk}\ ,
\eea
where we have introduced diagonal matrices $X,W$ 
with the nodes $x_i$ and the weights $w_i$ 
along their diagonals, and the matrix $Q$ with 
entries $Q_{ik}=q_i(x_k)$. More concisely, we have 
\be  B = QWQ^T \ , \qquad  A = QWXQ^T \ . \la{meqs}\ee 
Inserting these into the second relation in (\r{VD}) we obtain
\be QWXQ^T V = QWQ^T V D\ , \ee
and therefore
\be X = (Q^T V) D (Q^T V)^{-1}\ . \ee
Since diagonal matrices are conjugate if and only if they are equal,
up to reordering of the diagonal entries, we
deduce the first statement in the theorem. Furthermore, after
reordering the nodes if necessary, we see that the 
matrix $Q^TV$ commutes with $X$. But $X$ is a diagonal matrix with 
distinct entries on the diagonal, so $Q^TV$ must also be diagonal. 
Write $\sqrt{W}Q^TV=Y$. Then $Y^TY= V^T Q W Q^T V=I$ (by (\r{meqs}) and the first 
relation in (\r{VD})). So the diagonal entries of $Y$ are all plus or 
minus $1$. Writing $\sqrt{W}Q^T=YV^{-1}$ we obtain 
$w_iQ_{ji}^2=((V^{-1})_{ij})^2$ for all $i,j$, giving the 
second statement in the theorem. 
$\bullet$

\noindent{\bf Notes:}
1. Under a change of the basis $\{q_i\}$ we will have 
$Q\rightarrow MQ$, $A\rightarrow MAM^T$, $B\rightarrow MBM^T$, 
$V\rightarrow M^{-T}V$,  where $M$ 
is a suitable nonsingular matrix. The nodes and the weights,
 however, remain unchanged. 

\noindent 2.
Once the nodes $x_i$ are known, then the 
weights can be obtained from the 
first equation in (\r{meqs}). 
This, however, requires full knowlegde of the basis $q_i$,
whereas the formula for the weights in the theorem requires
knowledge of just a single basis element.

\noindent 3.
The simplest cases of theorem 3 are well-known. If we choose
$q_i(x)$ to be the orthormal basis of ${\cal P}_n$ formed by
Gram-Schmidt orthonormalization of the basis $1,x,\ldots,x^n$,
then $B=I$ and $A$ is tridiagonal. The method for finding nodes
and weights reduces to the classic and widely used method of [6]
(see [1] for an exposition). 
If we choose $q_i(x)=x^{i-1}$ then the matrices $A$ and $B$ are 
Hankel matrices ($A_{ij}$ and $B_{ij}$ both depend only on $i+j$),
and the computation of $D$ and $V$ is numerically unstable (at least
for all but the smallest values of $n$). For more general bases 
the method is essentially equivalent to that of [7] and [5], but the 
formulation here makes the issue of basis independence much 
clearer. 
The papers [7] and [5] focus on bases in which the 
$q_i(x)$ satisfy a recursion of the form
\be xq_i(x) = a_iq_{i+1}(x) + b_iq_i(x) + c_iq_{i-1}(x)\ . \ee
This makes most of the entries of $A$ expressible in terms of the 
entries of $B$ and the coefficients $a_i,b_i,c_i$. 

\noindent{\bf Example:} For the weight function $w(x)=(1+x)^{-1}$ on
$[0,1]$, the construction of orthogonal polynomials by Gram-Schmidt 
orthonormalization of the standard basis $1,x,\ldots,x^n$ is awkward.
It is much more convenient to work in  a basis in which most
of the elements have a factor $(1+x)$. As a first attempt, take
\be q_i(x) = \left\{
   \matrix{ (1+x)x^{i-1}  &  i=1,\ldots,n \cr
            1         &   i=n+1   \cr} \right. \ . \ee
It is straightforward to compute the necessary integrals to obtain
\be B_{ij} = \left\{ 
   \matrix{ \frac1{i+j-1} + \frac1{i+j}   &  i,j=1,\ldots,n \cr
             \frac1{i}  &  i=1,\ldots,n\ ,j=n+1 \cr
             \frac1{j}  &  j=1,\ldots,n\ ,i=n+1 \cr
            \ln 2     &   i=j=n+1 \cr }  \right. \quad , \ee
\be A_{ij} = \left\{ 
   \matrix{ \frac1{i+j} + \frac1{i+j+1}   &  i,j=1,\ldots,n \cr
             \frac1{i+1}  &  i=1,\ldots,n\ ,j=n+1 \cr
             \frac1{j+1}  &  j=1,\ldots,n\ ,i=n+1 \cr
            1-\ln 2     &   i=j=n+1 \cr }  \right.  \quad . \ee
For fairly small $n$, it is straightforward to 
find $D,V$ and the nodes and weights using Matlab. 
As $n$ gets larger and the conditioning of the Hankel 
matrices $B$ and $A$ becomes 
a concern, it is preferable to use the basis
\be q_i(x) = \left\{
   \matrix{ (1+x)p_i(x)  &  i=1,\ldots,n \cr
            1         &   i=n+1   \cr} \right. \ee
where the $\{p_i(x)\}$ are the  standard orthonormal polynomials on $[0,1]$
(orthonormal with respect to the weight function $1$).
Then $B$  and $A$ are tridiagonal and pentadiagonal matrices respectively. 

The point in this example is that we can exploit the freedom 
of basis to simplify  evaluation of the necesary integrals. 
The cost is a little more linear algebra, but this is hardly a challenge 
by modern standards. 

\noindent {\bf Generalization 1. Gaussian Quadrature on the Unit Circle} [2]. 
Let $w(\theta)$ be a  $2\pi-$periodic non-negative weight function. 
Assume there exists a quadrature formula 
\be  
\frac1{2\pi}\int_0^{2\pi} f\left(e^{i\theta}\right)  w(\theta) d\theta \approx 
\sum_{r=1}^{n+1} w_rf(z_r) \ ,  \qquad w_r\not=0\ ,
\label{quad3}
\ee
with distinct nodes $z_r$, 
which is exact for the $2n+2$ functions $f(z)=z^{-n},\ldots,z^n,z^{n+1}$ 
and all their linear combinations.
Let $q_r(z)$, $r=1,\ldots,{n+1}$ be a basis for the vector space of  (complex)
polynomials of degree at most $n$ in $z$. 
Define the $(n+1)\times(n+1)$ matrices $B,A$  by
\bea
B_{rs}&=&\frac1{2\pi}\int_0^{2\pi}  w(\theta)  
     \overline{q_r(e^{i\theta})}  q_s(e^{i\theta})  
      d\theta \ ,\qquad r,s=1,\ldots,{n+1}\ , \\
A_{rs}&=&\frac1{2\pi}\int_0^{2\pi}  w(\theta)  
     \overline{q_r(e^{i\theta})}  e^{i\theta} q_s(e^{i\theta})  
      d\theta \ ,\qquad r,s=1,\ldots,{n+1}\ .
\eea
Note that $B$ is Hermitian and positive definite, but $A$ need not be Hermitian. 
We assume the pair $A,B$ is diagonalizable,
i.e. that there exists a diagonal matrix $D$ and a matrix $V$ such that
\be 
A V = B V D \ . \la{VD2}
\ee
The matrices $D,V$ can be obtained in Matlab by the command {\tt [V D]=eig(A,B)}.
We then have the following analog of  theorem 3:

\noindent{\bf Theorem 3$'$}.  
The nodes and weights of the quadrature formula (\r{quad3}) (that is exact
for $f(z)=z^{-n},\ldots,z^n,z^{n+1}$) are determined by
\be
z_r = D_{rr}\ , \qquad
w_r = \frac{(V^{-1})_{rs}(BV)_{sr} }
           {q_s(z_r)\overline{q_s\left(\frac1{\overline{z_r}}\right)}}  
    \ . \ee
(The second relation holds for arbitrary $s$). 

\noindent{\bf Proof}: Proceeding as in the real case we obtain 
\be B=QW\tilde{Q}\ , \qquad A=QWX\tilde{Q} \ ,  \la{bqw}\ee
where $Q,\tilde{Q}$ are matrices with entries $Q_{rs}=q_r(z_s)$ and 
$\tilde{Q}_{rs}=\overline{q_s\left(\frac1{\overline{z_r}}\right)}$, 
respectively, 
and $X,W$ are diagonal matrices with diagonal entries equal to the nodes 
and the weights, respectively.
Using these relations in (\r{VD2}) we rapidly obtain
\be X=(\tilde{Q}V) D (\tilde{Q} V)^{-1}\ , \ee
giving us the first result, that $X=D$, and also that 
the matrix $Y=\tilde{Q}V$ is diagonal. 

From $\tilde{Q}=YV^{-1}$ we deduce that
\be  \tilde{Q}_{rs}=Y_{rr}(V^{-1})_{rs} \ , \la{ty2} \ee
for each $r,s$. 
Multiplying the first equation in (\r{bqw}) on the right by $V$, gives
$BV=QWY$, implying 
\be  (BV)_{sr}=Y_{rr}w_r Q_{sr} \ ,\la{ty1}  \ee
for each $r,s$. Finally we divide (\r{ty1}) by (\r{ty2}) to obtain the 
second result in the theorem.  $\bullet$ 

\noindent{\bf Notes}: 1. The result in the theorem shows how to compute
the weights given $A,B$ and one of the $q_r$. It may be easier, if
all of the $q_r$ are known, to directly use the first formula in (\r{bqw}) to
find the weights. 

\noindent 2. For each $k$ we have $A{\bf v}_k = z_kB{\bf v}_k$, 
where ${\bf v}_k$ denotes the 
vector that is the $k$th column of $V$. Thus ${\bf v}_k^*(A-z_kB){\bf v}_k=0$, or
\be \frac1{2\pi}\int_0^{2\pi} (e^{i\theta}-z_k) |a(\theta)|^2 
   w(\theta) d\theta = 0\ , \ee
where $a(\theta)=\sum_r ({\bf v}_k)_r q_r(e^{i\theta})$. Since
\be \left| \int_0^{2\pi} e^{\sqrt{-1}\theta}|a(\theta)|^2 
   w(\theta) d\theta \right|
\le 
\int_0^{2\pi}|a(\theta)|^2 
   w(\theta) d\theta \ ,
\ee
it follows that $|z_k|\le 1$, i.e. the nodes are inside the unit circle. 

\noindent{\bf Example:} Taking $w(\theta)=\sin^2\theta$ and using 
the standard basis $\{q_r(z)\}=\{1,z,\ldots,z^n\}$ we have 
\bea
A_{rs}&=&\frac1{2\pi}\int_0^{2\pi} 
   e^{i(r-s)\theta} \sin^2\theta
      d\theta =\frac14\left(2\delta_{r-s}-\delta_{r-s+2}-\delta_{r-s-2}\right)
   \ ,\\
B_{rs}&=&\frac1{2\pi}\int_0^{2\pi} 
   e^{i(r-s+1)\theta} \sin^2\theta      d\theta
=\frac14\left(2\delta_{r-s+1}-\delta_{r-s+3}-\delta_{r-s-1}\right)
  \ ,
\eea
where $\delta_r$ is $1$ if $r=0$ and $0$ otherwise, and we have 
used $\frac1{2\pi}\int_0^{2\pi} e^{i r\theta} d\theta
=\delta_r$. With $n=7$ and these 
simple choices of $A$ and $B$, and with a little help from 
Matlab, theorem 3$'$ rapidly
reproduces the nodes and weights given in the 
third part of table 1 in [2]. 

\noindent
For a further reference on the use of Gaussian quadrature for 
complex integrals see [8]. 

\noindent{\bf Generalization 2. Gauss-type rules  with some fixed
nodes.} 
We seek nodes and weights $x_i,w_i$, $i=1,\ldots,n+1$, and 
weights $v_\alpha$, $\alpha=1,\ldots,m$, such that the quadrature formula
\be
\int_a^b f(x) w(x) dx \approx \sum_{\alpha=1}^m v_\alpha f(y_\alpha) +
\sum_{i=1}^{n+1} w_if(x_i) \ , \qquad w_i,v_\alpha\not=0\ ,
\label{quad2}
\ee
is exact when $f(x)$ is any polynomial of degree not more than $2n+m+1$. 
Here the $y_\alpha$, $\alpha=1,\ldots,m$,  are given nodes 
in $[a,b]$, and we assume the $x_i$ and $y_\alpha$  are all distinct.
There is no particular reason to assume $w(x)$ is non-negative. 
Let $q_i(x)$, $i=1,\ldots,n+1$ be an arbitrary
basis for the vector space of polynomials of degree at most $n$. Let
the $(n+1)\times(n+1)$ matrices $B,A$ be defined by
\bea
B_{ij}&=&\int_a^b w(x) \prod_{\alpha=1}^m (x-y_\alpha)\ 
 q_i(x) q_j(x) dx\ ,\qquad i,j=1,\ldots,n+1\ , \\
A_{ij}&=&\int_a^b w(x) \prod_{\alpha=1}^m (x-y_\alpha)\ 
          x q_i(x) q_j(x) dx\ ,\qquad i,j=1,\ldots,n+1\ .
\eea
Both $B$ and $A$ are symmetric, but in this generalization $B$ need not
be positive  definite (though in the important cases
$m=1$, $y_1=a$ or $b$, and $m=2$, $y_1=a$, $y_2=b$,
it is positive or negative definite). 
Assume the pair $A,B$ is diagonalizable, i.e. that we can find 
a diagonal matrix $D$ and a matrix $V$ satisfying
\be 
A V = B V D \ . \la{VD3}
\ee
We then have the following: 

\noindent{\bf Theorem 3$''$}.  
If the quadrature formula (\r{quad2}) is exact
when $f(x)$ is any polynomial of degree not more than $2n+m+1$,
and the pair $A,B$ is diagonalizable, then 
the nodes and weights are determined by
\be
x_i = D_{ii}\ , \qquad
w_i\prod_{\alpha=1}^m (x_i-y_\alpha)
 = \frac{(V^{-1})_{ij}(BV)_{ji} }
           {q_j(x_i)^2}  
    \ . \ee
(The second relation holds for arbitrary $j$). 

The proof of this is sufficiently similar to that of theorem 3$'$ that we omit
it. 

Once the nodes $x_i$ 
and weights $w_i$ have been determined using this theorem, the 
weights $v_\alpha$ can be determined by
solving the linear system arising from exactness of the 
quadrature formula for $f(x)=1,x,\ldots,x^{m-1}$. 

\vskip.4in
\noindent{\bf Acknowledgement:} We thank David Tannor for repeatedly 
pushing us to understand the results in [4], which led to this 
work. We thank Gene Golub and Walter Gautschi for useful comments.
This letter is a revised version, with additions, of a manuscript with the same title
by the second author written in April 2003 but not submitted for 
publication. 

\medskip
\noindent{\bf References.}

\noindent[1] P.J.Davis and P.Rabinowitz, {\em Methods of Numerical Integration},
Academic Press (1975).

\noindent[2] L.Darius, P.Gonz\'alez-Vera and F.Marcellan,
{\em Gaussian Quadrature Formulae on the Unit Circle},
{\sl J.Comp.Appl.Math.} {\bf 140} (2002) 159-183.

\noindent[3] I.Degani, J.Schiff and D.J.Tannor,
{\em Commuting Extensions and Cubature Formulae},
{\sl Numer.Math.}, to appear. arXiv:math.NA/0408076.

\noindent[4] A.S.Dickinson and P.R.Certain, {\em Calculation of 
Matrix Elements for One-Dimensional Quantum-Mechanical Problems},
{\sl J.Chem.Phys.} {\bf 49} (1968) 4209-4211.

\noindent[5] W.Gautschi, {\em On the Construction of Gaussian Quadrature Rules 
from Modified Moments}, {\sl Math.Comp.} {\bf 24} (1970) 245-260. 

\noindent[6] G.H.Golub and J.H.Welsch, {\em Calculation of 
Gauss Quadrature Rules}, {\sl Math.Comp.} {\bf 23} (1969)
221-230. 

\noindent[7] R.A.Sack and A.F.Donovan, {\em An Algorithm for 
Gaussian Quadrature 
given Modified Moments}, {\sl Numer.Math.} {\bf 18} (1972) 465-478.

\noindent[8] P.E.Saylor and D.C.Smolarski, {\em Why Gaussian 
Quadrature in the Complex Plane?}, {\sl Numerical Algorithms}
{\bf 26} (2001) 251-280. 

\end{document}